\documentclass [11pt] {amsart}
\usepackage{amsthm,amsmath,amssymb,amscd,graphics,enumerate}
\usepackage{latexsym}
\usepackage{verbatim}
\usepackage{enumerate}
\usepackage{amsthm}
\usepackage{amsmath}
\usepackage{amscd}

\newcounter{theorem}[section]
\numberwithin{equation}{section}

\newtheorem{cl}[theorem]{Claim}
\newtheorem{subcl}[theorem]{Subclaim}

\newtheorem{Thm}[theorem]{Theorem}
{\theoremstyle{remark}
 
\newtheorem{Rem}[theorem]{\text{\textbf{Remark}}} }
\newtheorem{Def}[theorem]{Definition}
\newtheorem{Lem}[theorem]{Lemma}
\newtheorem{Prop}[theorem]{Proposition}
\newtheorem{Cor}[theorem]{Corollary}

\newtheorem{notation}[theorem]{Notation}
\newtheorem{assumptions}[theorem]{Assumptions}
\newtheorem{Ex}[theorem]{Example}
{\theoremstyle{definition}
 }

\newcommand{\CC}{{\mathbb{C}}}

\newcommand{\PP}{{\mathbb{P}}}

\newcommand{\mO}{\mathcal{O}}
\newcommand{\mX}{\mathcal{X}}

\newcommand{\mC}{\mathcal{C}}

\newcommand{\mD}{\mathcal{D}}
\newcommand{\mE}{\mathcal{E}}

\newcommand{\mH}{\mathcal{H}}

\newcommand{\mK}{\mathcal{K}}

\newcommand{\mU}{\mathcal{U}}

\newcommand{\mZ}{\mathcal{Z}}
\newcommand{\fbar}{\bar{f}}

\title{Characterizing  projective spaces for varieties with at most quotient singularities}
\author{Jiun-Cheng Chen}
\address{Department of Mathematics\\ Northwestern
University\\ 2033 Sheridan Road\\ Evanston\\ IL 60208-2370\\ USA}
\email{jcchen@math.northwestern.edu}

\date{\today}

\begin{document}
\begin{abstract}
We generalize the well-known numerical criterion for projective spaces  by Cho, Miyaoka and Shepherd-Barron  to  varieties with at worst quotient singularities.    
Let $X$ be a  normal projective variety of dimension $n \geq 3$ with at most quotient singularities.  Our result asserts that if $C \cdot (-K_X) \geq n+1$ for every curve $C \subset X$, then $X \cong \PP^n$. 
\end{abstract}
\maketitle
\markboth{}{}
\section{Introduction}
We work over the field $\CC$. 
The $n$-dimensional  projective space $\PP^n$ is  probably the simplest  compact (projective)   complex manifold. 
Let $K_{\PP^n}$ be the canonical bundle (the line bundle of holomorphic $n$-forms).  
It is elementary to see that the  $-K_{\PP^n}=\mO(n+1)$. In particular,  $(line) \cdot (-K_X)=n+1$    and $C \cdot (-K_{\PP^n}) \geq n+1$ for all curves $C \subset \PP^n$.   
This is  a rather  unusual property:   recall that
 if  $X$ is a smooth  projective variety of dimension $n$ and $K_X$ is 
not nef,   
 then  the Cone theorem  implies that 
  there is a rational curve $C$ such that $ 0<C \cdot (-K_X) \leq n+1$.

From this perspective, the anti-canonical bundle $-K_{\PP^n}$ is unusually positive.  
It turns out that  this property characterizes $\PP^n$ among smooth projective varieties of dimension $n$.
\begin{Thm}[\cite{cmsb02} and \cite{ke01}]\label{cmsb}
Let $X$  be a smooth  projective variety of dimension $n \geq 3$. Assume that
$C \cdot (-K_X) \geq n+1$ for all curves  $C \subset X$.
Then $X\cong \PP^n$.
\end{Thm}

This result was first proved by Cho, Miyaoka and Shepherd-Barron \cite{cmsb02} and later by Kebekus \cite{ke01}.   

Since the condition $$ C \cdot (-K_X) \geq n+1 \; \forall \; C \subset X$$ is very strong, it is natural to ask if the  assumption on smoothness  is necessary. 
In \cite{ct05a}, H,-H, Tseng and the author proved the characterization result  assuming that  $X$ has only isolated LCIQ singularities.   The main goal of this paper is to further weaken  the smoothness assumption;  it  
 builds upon methods developed in \cite{ct05a}. 
The precise statement of our result is the following:
\begin{Thm}\label{main} Let $X$ be a projective variety
of dimension $n \geq 3$ with at most quotient  singularities.
Assume that $$C \cdot  (-K_X) \geq n+1 \; \forall \; C  \subset X.$$ 
Then $X \cong\PP^n$.
\end{Thm}

  We now explain the strategy used in \cite{ke01}. 
First consider the projective space $\PP^n$.  Let $p \in \PP^n$ be any point. 
Let $\widetilde{\PP^n}$ be the blowing up  of $\PP^n$ along the point $p$.
The  variety $\widetilde{\PP^n}$  is  a $\PP^1$-bundle over $\PP^{n-1}$.  
Let $E \cong \PP^{n-1} \subset \widetilde{\PP^n}$ be the exceptional divisor. 
The normal bundle of the exceptional divisor $E \cong \PP^{n-1}$ is
$\mO_{\PP^{n-1}}(-1)$.  The variety $\widetilde{\PP^n}$  is   the 
 Chow family of 
minimal degree  rational curves  (lines)  through the point $p \in \PP^n$.

Now consider the variety $X$. 
Take a  general point $x \in X$. Denote by $\tilde{X}$ the blowing up of $X$
along $x$.  
 If we can prove that $\tilde{X}$ is the Chow family of  minimal degree (with respect to an ample line bundle)  rational curves through $x$,  then we have a good chance to prove that  $\tilde{X} \cong \widetilde{ \PP^n}$ and hence
$X \cong \PP^n$.
  
When $X$ is smooth, 
Kebekus  proved that  the Chow scheme $H_x$ of  minimal degree rational curves through  a general point $x \in X$
is isomorphic to $\PP^{n-1}$ and $\tilde{X} \cong \widetilde{ \PP^n}$
 \cite{ke01}. It follows easily that $X \cong \PP^n$. 
One important  step  in  his proof   is to show that  $X$ has a lot of
minimal degree rational curves through $x$. More precisely, one needs to show
that 
$$dim\;H_x \geq  l \cdot (-K_X) -2=n-1,$$ 
where $[l] \in H_x$.
Note that $dim\; H_x \leq n-1$ by a standard bend and break argument.  It follows that   the Chow scheme  $H_x$ 
has the expected dimension $n-1$.   
Kebekus  then  proved that the tangent map 
$H_x \to \PP^{n-1}$ is  an isomorphism \cite{ke00} \cite{ke01}.     

When $X$ is possibly singular, the situation is quite different: it is
difficult to have the desired lower bound on the 
dimension of $H_x$.   
In \cite{ct05a} (joint work with H.-H. Tseng), we  proved  the
characterization result  when $X$ has at 
worst isolated LCIQ  singularities.
We used   the Deligne-Mumford stack $\mX \to X$, twisted stable maps into $\mX$  and Mori's bend and break
techniques to show the existence of 
a rational curve $f: \PP^1 \to X$ such that (1) 
$f(\PP^1) \cdot (-K_X) =n+1$, and (2) $f(\PP^1)$ does not meet the singular locus
of $X$. 
Using this fact,  one can show  (following \cite{ke01}) that  $H_x \cong \PP^{n-1}$ and $X \cong \PP^n$.   
In \cite{ct05a}, 
 the condition that $X_{sing}$ is  isolated  is essential; it ensures that we
 can deform a specific rational curve  in $X$.  
 Our methods do not 
 apply without  this assumption.

In this paper, we develop  methods  to study  the case when $X$ has only quotient singularities (not necessarily isolated). 
The main idea is quite simple: instead of considering the twisted stable map directly, we consider a double cover and  study  the possible degeneration types. 
 The assumption $$C \cdot (-K_X) \geq n+1\; \forall \; C \subset X $$
 ensures that the possible degeneration types are  limited and 
we are able to show the existence of a minimal degree rational curve $f:\PP^1 \to X$ which does not meet $X_{sing}$. 
Once this fact  is established, it is  standard \cite{ke01} \cite{ct05a} to  prove that the Chow family has the right dimension, i.e. $n-1$,  and  $X \cong \PP^n$.

A few words on using Deligne-Mumford stacks: for a  singular variety $X$, it 
is difficult to  prove the existence of  enough   rational curves via
deformation theory; obtaining a 
lower bound of the dimension $Mor(\PP^1, X)$ (in terms of $-K_X$-degree) is difficult.   
However, if $X$ has only local complete intersection quotient (LCIQ)  singularities, 
 one can study representable morphisms from a {\em twisted curve}\footnote{Roughly speaking, this is a one dimensional Deligne-Mumford stack with isolated cyclic quotient stack structures and with coarse moduli space being a nodal curve \cite{av02}.} to the stack $\mX$ whose coarse moduli space is $X$. This provides  a reasonable alternative. 
We can obtain a lower bound on the dimension of  the space of twisted stable
maps expressed 
in terms of the $-K_X$-degree and the number of twisted points, see \cite{ct05}.

The rest of this paper is organized as follows:
In Section 2, we recall basic definitions on  twisted curves and twisted stable maps. We also give a formula on the lower bound on the morphism space.  The main proposition (Proposition ~\ref{general}) is proved in Section 3. 
In Section 4, we
 sketch the proof following \cite{ke01};  we do not claim any originality of these results.  
 We also make a few  remarks in 
 that section.
\section*{Acknowledgments}  Part of this research was conducted  while 
the author   was attending  the JAMI 
conference at Johns Hopkins University. He   likes to
thank the conference organizers.  
He also likes to thank Dan Abramovich,  Lawrence Ein and Stefan Kebekus 
for helpful 
discussions and valuable suggestions.  
\section{Twisted  curves and twisted stable maps}\label{twisted}
Twisted curves play an important role in this paper. One of the main motivations of introducing twisted curves
 is    to compactify the space (Deligne-Mumford stack) of stable maps into a
  proper Deligne-Mumford stack \cite{av02}.   
Roughly speaking, twisted curves  are nodal curves having  certain stack structures \'etale locally near nodes (and, for pointed curves, marked points).  
For the precise definition, see
 \cite{av02}, Definition 4.1.2.

Let $\mC$ be a twisted curve and $C$ its coarse moduli space.
\subsubsection{Nodes}
For a positive integer $r$, let $\mu_r$ denote the cyclic group of $r$-th roots of unity. \'Etale locally near a node, a twisted curve $\mC$ is isomorphic to the stack quotient $[U/\mu_r]$ of the nodal curve $U=\{xy=f(t)\}$ by the following action of $\mu_r$: $$(x,y)\mapsto (\zeta_r x,\zeta_r^{-1} y),$$ where $\zeta_r$ is a primitive $r$-th root of unity. \'Etale locally near this node, the coarse curve $C$ is isomorphic to the schematic quotient $U/\mu_r$.
\subsubsection{Markings}
\'Etale locally near a marked point, $\mC$ is isomorphic to the stack quotient $[U/\mu_r]$. Here $U$ is a smooth curve with local coordinate $z$ defining the marked point, and the $\mu_r$-action is defined by $$z\mapsto \zeta_r z.$$ Near this marked point the coarse curve is the schematic quotient $U/\mu_r$.
\subsection{Twisted stable maps}
\begin{Def}
A twisted $n$-pointed stable map of genus $g$ and degree $d$ over a scheme $S$ consists of the following data (see \cite{av02}, Definition 4.3.1):
$$\begin{CD}
\mC @>f>>  \mX \\
@V{\pi_\mC}VV @V{\pi}VV \\
C @>\fbar>> X \\
@V{}VV \\
S.
\end{CD}$$
along with $n$ closed substacks $\Sigma_i\subset \mC$ such that
\begin{enumerate}
\item $\mC$ is a twisted nodal $n$-pointed curve over $S$ (see \cite{av02}, Definition 4.1.2),
\item $f:\mC\to \mX$ is representable,
\item $\Sigma_i$ is an \'etale gerbe over $S$, for $i=1,...,n$, and
\item the map $\fbar: (C,\{p_i\})\to X$ between coarse moduli spaces induced from $f$ is a stable $n$-pointed map of degree $d$ in the usual sense.
\end{enumerate}
\end{Def}
A twisted map  $f:\mC\to \mX$ is stable if and only if for every irreducible component $\mC_i\subset \mC$, one of the following cases holds:
\begin{enumerate}
\item
$f|_{\mC_i}$ is nonconstant,
\item
$f|_{\mC_i}$ is constant, and $\mC_i$ is of genus at least $2$,
\item
$f|_{\mC_i}$ is constant, $\mC_i$ is of genus $1$, and there is at least one special points on $\mC_i$,
\item
$f|_{\mC_i}$ is constant, $\mC_i$ is of genus $0$, and there are at least three special points on $\mC_i$.
\end{enumerate}
In particular, a nonconstant representable morphism from a smooth twisted curve to $\mX$ is stable.

We say a twisted stable map $\mC \to \mX$ is rational if the coarse moduli space $C$ of $\mC$ is rational.

Let $\mK_{g,n}(\mX,d)$ denote the category of twisted $n$-pointed stable maps to $\mX$ of genus $g$ and degree $d$. The main result of \cite{av02} is that $\mK_{g,n}(\mX,d)$ is a proper Deligne-Mumford stack with projective coarse moduli space denoted by  $K_{g,n}(\mX,d)$.
Let $\beta \in H_2(X)$ be a homology class. 
The space of twisted $n$-pointed stable maps $f:\mC \to \mX$ of genus $g$  and  
  homology class $[(\pi \circ f)_*(\mC)]= \beta$ is denoted by  $\mK_{g,n}(\mX, \beta)$. 
This  stack is also proper \cite{av02}.
\subsection{Morphism space  from a twisted curve to a Deligne-Mumfors stack}
In this paper, we use both the stack of  twisted stable maps and  
the morphism space from $\mC$ to $\mX$. Roughly speaking, an element in 
the morphism space $Mor(\mC, \mX)$ is a twisted stable map together with a 
parameterization
on the source curve $\mC$.   
Let  $\Sigma \subset \mC$ be the set of twisted points and 
$B \subset \mC$  a finite set of points (twisted or untwisted).
Let $f:\mC \to \mX$ be a representable morphism.  
When $\mX$ is smooth, 
we have a lower bound on  
the dimension of  $Mor(\mC,\mX;f|_B)$ near the morphism $[f]$.
\begin{Lem}[= \cite{ct05} Lemma 4.4] \label{def}
$$dim_{[f]}Mor(\mC,\mX;f|_B)\geq -\mC\cdot K_\mX + n [\chi(\mO_{\mC}) - \;Card(B)]-\sum_{x\in \Sigma\setminus B} age(f^*T\mX,x).$$
\end{Lem}
\begin{Rem}
When $\mX$ has only  LCI singularities, a
similar formula still holds as long as the image of $\mC$ does not lie completely in  $\mX_{sing}$ (the locus where  the  stack $\mX$ is singular).     
\end{Rem}

\subsection{Lifting}
Let $X$ be a normal projective variety with quotient singularities 
and $\mX$ a  proper smooth Deligne-Mumford 
stack such that  $\pi: \mX \to X$ is 
isomorphic over $X_{reg}=X-  X_{sing}$ and 
$X$ is a coarse moduli space of $\mX$.
 Let $C$ be a smooth irreducible curve  and $\fbar: C \rightarrow X$ a morphism. 
 
We want to ``lift''  the map $\fbar:C \to X$ to a map
$C \to \mX$.  In general, this is not possible unless we endow a orbicurve 
structure on  $C$ \cite{av02} Lemma 7.2.5,  or pass to  a  finite cover  $C' \to C$ \cite{av02} Theorem 7.1.1.

First consider the case when the image  $\fbar(C)$ meets the smooth locus of $X$. 
 Let $\{p_i| i \in I\}\subset C$ be the finite set of points which are mapped to the
singular locus of $X$, and let $C_0=C\setminus \{p_i| i \in I\}$. Since $X$ is 
isomorphic to $\mX$ away from  the singular locus $X_{sing}$, the map $\fbar|_{C_0}:C_0\to X$ admits a lifting $C_0\to \mX$. By \cite{av02}, Lemma 7.2.5, there exists a twisted curve $\mC$ with coarse moduli space $C$, and a twisted stable map $f:\mC\to \mX$ extending $C_0\to \mX$. 

Now consider the case when  $\fbar(C) \subset X_{sing}$.  
Let $\eta \in C$  be the generic point of $C$. 
Consider the morphism
$\eta  \in C \to X$.  After a finite  extension $\eta' \to \eta$,  
there is a lifting  $\eta ' \to \mX$ since $\mX$ is proper. 
Let $\eta ' \in C'$ be the smooth irreducible curve.  Note that  the morphism $\eta ' \to \mX$ is defined over an open set  $U \to \mX$.
One can extend the finite morphism $\eta ' \to \eta$ to $C' \to C$ since $C$ is proper.
Let  $C' -U=
\{ p_j| j \in J\}$.  
 Endowing stack structures on this finite set of points, we can 
extend the morphism $U \to \mX$ to  $ \mC' \to \mX$. After a further 
finite cover  $C'' \to C'$, we can find a morphism
$C'' \to \mX$.

\section{Untwisted rational curves}
In this section, 
we prove the existence of  "good" rational curves with  $-K_X$-degree $n+1$, i.e. rational curves which do not intersect the singular locus of $X$.  

Consider a smooth Fano   variety $X$.   
Let $x \in X$ be a general point and  $S=\{x_1, x_2, \cdots, x_k \} \subset X$
be any finite set of points which does not contain $x$. 
It is  possible to  find a   rational curve through the  point $x$ which  misses the finite set $S=\{x_1, x_2, \cdots, x_k \} \subset X$. This property  is, however, not true for singular varieties.

\begin{Ex}
Let $E \subset \PP^2$ be a smooth elliptic curve. 
Let $X \subset \PP^3$ be  the projective cone of $E$. 
The surface $X$ has only  one 
LCI singularity at the vertex.  The surface $X$ is Fano and 
all rational curves pass through the vertex.  
\end{Ex}
This example shows that the existence of a rational curve which does not meet the singular locus  $X_{sing}$ is  non-trivial.

\begin{notation}\label{stack}
Let $X$ be a normal projective variety with at worst quotient singularities. 
Fix  a  proper smooth Deligne-Mumford stack $\pi: \mX \to X$ such that $X$ is a coarse moduli space of $\mX$ and  $\pi$ is 
an isomorphism over $X_{reg}=X \setminus X_{sing}$.  
Note that  $K_{\mX}=\pi^{*}K_X$ and $$\mC \cdot K_{\mX}= C
\cdot K_X$$ for any (twisted) curve $\mC \to \mX$ with coarse curve $C$
\cite{ct05}.
\end{notation}

 We start with the following  lemma: 
\begin{Lem}\label{1twisted}
Let $R \subset \overline{NE}(X)$ be any  $K_X$-negative  extremal ray. Then there exists a twisted rational
 curve $f:\mC \to \mX$ 
such that (1) $\mC$ has  at most one twisted point, (2) the intersection number  $\mC \cdot (-K_{\mX}) \leq n+1$, and (3) $[(\pi \circ f)_*(\mC)] \in R$.
\end{Lem}
 \begin{proof}
This is essentially  Proposition 3.1 in \cite{ct05a}; one only needs to note that  the argument   goes through when  
the stack $\mX$ has only isolated LCI singularities.   
\end{proof}

Unless mentioned  otherwise,
we make further assumptions in the rest of this paper: 
\begin{assumptions}\label{assu} From now on we assume that $X$ is a
projective variety of $dim X=n \geq 3$ with at most quotient 
singularities and has the property that 
$$C \cdot (-K_X) \geq n+1$$ for
 every curve $C \subset X$.
\end{assumptions}

\begin{Lem}\label{twisted1general}
Assumptions as in Assumptions \ref{assu}. Let  $R \in \overline{NE}(X)$ be any $K_X$-negative extremal ray and
$x \in \mX$  a general point. Then there exists a twisted rational curve
$f:\mC \to \mX$ such that (1) $\mC$ has at most 1 twisted point, (2)  $x \in f(\mC)$, (3) $[(\pi \circ f)_*(\mC)] \in R \subset H_2(X)$,  and  (4) $\mC \cdot f^*(-K_{\mX})=n+1$. 
\end{Lem}
\begin{proof}
We only need to verify the condition (2) by Lemma ~\ref{1twisted}. 
Let $f: \mC \to \mX$ be a twisted rational 
curve as in Lemma ~\ref{1twisted} and  $\infty \in \mC$  the twisted point on $\mC$. 
By Lemma ~\ref{def},
$$dim_{[f]}\;Mor(\mC, \mX, f|_{\infty}, n+1) \geq (n+1)+(1-1)n=n+1.$$ 
Denote by  $\mH_{f(\infty)} \subset \mK_{0,1}(\mX, [(\pi \circ f)_* \mC])$   the stack of twisted stable maps 
through the 
twisted point $f(\infty)$. 

Abusing  the notation, we  also view    $[f] \in Mor (\mC, \mX,f|_{\infty}, n+1)$ as an element 
 of   $\mH_{f(\infty)}$.
Note that   $dim_{[f]}\;\mH_{f(\infty)} \geq n-1$.
  
Let $T \to \mH_{f(\infty)}$ be any quasi-finite morphism where  $T$ is    a  normal irreducible 
quasi-projective  
variety of dimension $n-1$  such that $[f] \in Im(T)$.
Consider 
the pull back family $U_T \to T$ of twisted stable maps and the 
 morphism  $i_T: U_T \to \mX$. Forgetting the stack  structures on the source curves and on  the target stack $\mX$, one obtains a family of stable maps (into $X$).   
Denote  this family 
of stable maps 
 by $\phi_T: \overline{U}_T \to T$ and the morphism (into $X$) by $\bar{i}_T: \overline{U}_T \to X$.

Note that the homology class $[(\pi \circ f)_*(\mC)]$ can not be written as
the sum of at least two curve classes by Assumptions ~\ref{assu}.
Let $q= (\pi \circ f) (\infty)$. 

Recall  that 
 $\mK_{0,1}(\mX, [(\pi \circ f)_*(\mC)])    \to  \mK_{0,1}(X, [(\pi \circ f)_*(\mC)])$ is quasi-finite \cite{av02}. 
Being the composition of two quasi-finite morphisms, the  morphism $T  \to
 \mH_{f(\infty)} \to \mK_{0,1}(X, [(\pi \circ f)_*(\mC)])$
is also  quasi-finite.     

\begin{cl}\label{chow}
The   morphism  $\bar{i}_T: \overline{U}_T \to  X$ 
is quasi-finite away from the preimage of $q$.
\end{cl}
\begin{proof}
This is a standard bend and break argument.
Note that every  fiber of the family $\phi_T: \overline{U}_T \to T$ is irreducible since the homology class 
 $[(\pi \circ f)_* \mC]$ is unbreakable. 
It follows that no fiber of $\phi_T: \overline{U}_T \to T$ is contracted by the morphism $\bar{i}_T$.

Suppose  that  $\bar{i}_T$ is not quasi-finite (away from the preimage of $q$). Then there is a 
curve $C \subset \overline{U}_T$ (may not be projective) such that 
 $\bar{i}_T(C)=q_1 \neq q$. 
Let   $\phi_T(C) \subset T$ be the image of $C$ and  $C' \to \phi_T(C)$ the normalization. 
Note that  the morphism  
$$C' \to \phi_T(C) \subset  T \to \mK_{0,1}(X, [(\pi \circ f)_*(\mC)])$$ 
is still quasi-finite. 
Pull back the family of stable maps  to the curve  $C'$.
Compactify this   family of stable maps; this is a 
 family (over a proper    curve $\bar{C'}$)  of stable 
 maps whose image contains $q$ and $q_1$.  
This contradicts  the unbreakable assumption  on the  homology  class $[(\pi \circ f)_*(\mC)]$  by Mori's bend and break.
 \end{proof}

By Claim ~\ref{chow}, it  follows that  the dimension of the 
image $\bar{i}_T (\overline{U}_T)$ (and the dimension  
of $(\pi \circ i_T)(U_T)$)
is $n$. 
Thus a general point $x \in \mX$ lies on the image of  a 
twisted  stable map  which satisfies  the conditions (1), (3), and (4).  
This concludes the proof.     
\end{proof}

\begin{Prop}\label{general} 
Assumptions as in Assumptions ~\ref{assu}.
Let $R \subset \overline{NE}(X)$ be any  $K_X$-negative extremal ray. There is a rational curve $f:\PP^1 \to \mX$ such that (1) the image $\pi \circ f(\PP^1)$ does not lie in $X_{sing}$, (2) $x \in f(\PP^1 )$ for a general point $x \in \mX$, 
(3) the  class $[\pi \circ f (\PP^1)] \in R$, and (4) $f^* \PP^1 \cdot (-K_{\mX})=n+1$.
\end{Prop}
\begin{proof}
By Lemma ~\ref{twisted1general}, there is  a twisted 
rational  curve $\mC \to \mX$ 
such  that (1) $\mC$ has 
at most one twisted point, (2) $x_0 \in f(\mC)$ ($x_0$ is a general point of $\mX$),  
  (3) $[(\pi \circ f)_* (\mC)] \in R \subset  H_2(X)$, and (4) $\mC \cdot f^*(-K_X)= n+1$. 
  
Let $\beta= [(\pi \circ f)_* \mC] \in H_2(X)$.  We fix this homology class in the rest of the proof. 

We may assume that $\mC$  does have a twisted point, 
denoted by $\infty \in \mC$.   Denote by $x_{\infty}=f(\infty)$  the image of $\infty$.  May assume that $f(0)=x_0$ 
 and   $x_0 \in \pi^{-1}(X_{reg})$. 
\newline \newline
Step 0. Double cover of $\mC$: 

Recall that  the coarse curve of $\mC$ is $\PP^1$. 
Choose a 2-to-1 cover $h: \PP^1 \to \PP^1$ such that $h(0)=0$, \;
$h(\infty)=\infty$  and $h$ is ramified at 
$\{0,\;\infty\} \subset \PP^1$, e. g. $h(z)=z^2$.   
Choosing a 
suitable stack 
 structure at $\infty \in \PP^1$,
 we can lift $h: \PP^1 \to \PP^1$ to $\mD \to \mC$ where $\mD$ is a twisted 
rational curve. We abuse the notation and 
still denote the morphism on the stack level by $h: \mD \to \mC$. 
Denote the stacky point on $\mD$ by $\infty$ and 
the preimage of $0 \in \mC$ on $\mD$ by $0$ (an untwisted point).
    
The composition $f \circ h: \mD \to \mX$ is a twisted stable map into $\mX$ with $-K_{\mX}$-degree $2n+2$. Note that the homology class 
 $[(\pi \circ f \circ h)_*(\mD)]=2 \beta \in H_2(X)$. 
\newline \newline
Step 1. Bend and break:

We first show that we can deform the curve $f \circ h: \mD \to \mX$. 
Let  $Mor (\mD, \mC, 2)$ be the space  of
of degree $2$ representable morphisms  from $\mD$ to $\mC$. Let 
$Mor (\mD, \mC, h |_{\{0,\; \infty \}}, 2) \subset Mor (\mD, \mC, 2) $ 
be the subspace consisted of  morphisms $h_1:\mD \to \mC$ such that $h_1(0)=0$ and $h_1(\infty)=\infty$.  
The space    
 $Mor (\mD, \mC, h |_{\{0,\; \infty \}},2)$  has 
 dimension $3$.

By Lemma ~\ref{def},  
$$dim_{[f \circ h]}\; Mor (\mC, \mX, (f\circ h)|_{\{0, \infty \}},2n+2)
\geq 2n+2+(1-2)n=n+2 > 3.$$ 

Thus  we can deform 
  the morphism $f \circ h: \mD \to \mX$ such that the image deforms  in $\mX$. 
This also implies  that there is a morphism $[f'] \in Mor( \mD , \mX, 2\beta)$ such 
 that $f': \mD \to \mX$ is birational to its image. 

Choose $[g] \in Mor( \mD , \mX, 2\beta)$ such that $g: \mD \to \mX$ is birational 
 to its image in $\mX$ and $g(0)=x_0$.   
Let $1 \in \mD$ be an untwisted point and $x_1=g(1)$ its image in $\mX$.

\begin{cl}\label{threepoints}
Replacing $[g] \in Mor( \mD , \mX, 2\beta)$ and  the point $x_1$  if necessary, we may 
assume that $x_1 \in \pi^{-1}(X_{reg})$,  $x_1 \neq x_0$,  and  $x_1$ does not lie on the image of $[h]$ where  
  $[h] \in \mK_{0,k}(\mX, \beta),\;k=1,\;2$ is any twisted stable map  through $x_0$ and $x_{\infty}$.  
\end{cl} 
\begin{proof}
Let  $\mK_1 \subset \mK_{0,1}(\mX, 2\beta)$ be the family of twisted stable maps whose image contains $x_0$ and $x_{\infty}$.
View the morphism  $g: \mD \to \mX$ as an element of $\mK_1$.
Note that 
$$dim_{[g]}\; \mK_1 = dim\; Mor_{[g]}\; (\mD, \mX, g|_{\{0,  \; \infty\}}, 2\beta) -dim \; Aut(\mD, \{0,\;\infty\})$$ 
$$ \geq 2n+2 -(1-2)n -1=n+1.$$
Here we use Lemma ~\ref{def} and the fact  $dim\;Aut(\mC, \{0,\; \infty\})=1$.  
The image of this family is at least a surface since
the image of $g: \mD \to \mX$ deforms.  

Consider all possible twisted stable maps $[h] \in \mK_{0,k}( \mX, \beta), \;k=1,\;2$ such that the points  $x_0$ and $x_{\infty}$ lie  on the image of $h$.  
By Lemma ~\ref{moribandb}, there are only finitely many such such twisted stable maps.
The image of all these twisted stable maps is one dimensional (in $\mX$).
The claim follows by choosing a suitable 
 $[g] \in \mK_1$.
%
\end{proof}
Let $\mK \subset   \mK_{0,1}(\mX, 2\beta)$ be  the substack  of twisted stable maps whose image 
contains $x_0$, $x_1$ and $x_{\infty}$.   
Compute the dimension of the stack  $\mK \subset   \mK_{0,1}(\mX, 2\beta)$ at $[g]$:  
$$dim_{[g]}\; \mK = dim\; Mor_{[g]}\; (\mD, \mX, g|_{\{0, \;1,\; \infty\}}, 2\beta) -dim \; Aut(\mD, \{0,\;1,\;\infty\})$$ 
$$ \geq 2n+2 -(1-3)n -0=2.$$
By Mori's bend and break,  
the domain curve of some twisted stable map  $[h] \in \mK$  has to degenerate
to  at least two irreducible components. 
\newline \newline Step 2.  Analysis on possible types of degeneration:
  
Since $\mD \cdot (-K_{\mX})= 2n+2$ and every curve 
has $-K_{\mX}$-degree at 
least $n+1$, only two components of the domain curve  are not contracted. 
Note that  the coarse space of domain curve is a rational tree with an extra special 
point (coming from the original 
stacky point of $\mD$). 
It is easy to see that the domain curve  can only break into two or three pieces. 
\newline \newline
Case I. The domain curve has three irreducible components. \newline \newline
Denote these components  by $\mD_1$, $\mD_2$ and $\mD_3$. 
Denote by $f_1: \mD_1 \to \mX$, $f_2: \mD_2 \to \mX$ and $f_3:\mD_3 \to \mX$ the twisted stable maps.
Let $\mD_2$ be the component which intersects other two components. Note that
$D_2$ has to be  contracted.  

Note that $[\pi \circ f_1 (\mD_1)]+ [\pi \circ f_3(\mD_3)]= [\pi \circ g (\mD_1)]=2 \beta$, and 
$$(\pi \circ f_1)_*\mD_1 \cdot (-K_{X})=n+1= (\pi \circ f_3)_*\mD_3 \cdot (-K_{X}).$$  
 Since $\beta \in R$ and $R$ is an extremal ray, the class  $[\pi \circ f_1 (\mD_1)]$ is a multiple of the 
class $[\pi \circ f_3(\mD_3)]$.  Since they have the same $-K_X$-degree,  $[\pi \circ f_1 (\mD_1)]= [\pi \circ f_3(\mD_3)]=\beta$.

By symmetry, 
we only need to consider the following  cases:
\newline \newline
Case I-a:  The stacky point $\infty \in \mD_2$ 
and $x_0 \in f_1(\mD_1)$ and $x_1 \in f_3( D_3)$.
\newline \newline
Case I-b: The stacky point $\infty \in \mD_2$ and 
$\{x_0,\;x_1\} \subset f_1(\mD_1)$.
\newline \newline
Case I-a.   We will show that  Case I-a forms a finite set.
Pushing forward to $X$, and 
noting  that 
\begin{enumerate}
\item  $\pi(x_0) \in (\pi \circ f_1)( \mD_1)$, \;  $\pi (x_{\infty})
\in (\pi \circ f_1)(\mD_1)$, and 
\item  
 $\pi(x_1) \in  (\pi\circ f_3)(\mD_3)$,\; $\pi(x_{\infty}) \in  \pi\circ
f_3(\mD_3)$,
\end{enumerate}
 it follows that  there are only finitely many such 
stable maps
by  
 Lemma ~\ref{moribandb}.   
\newline \newline
Case I-b. This  is impossible by our choice of $x_0,\;x_1$, and  $ x_{\infty}$
(see Claim ~\ref{threepoints}).
\newline \newline
Case II. The domain curve has two components, denoted by  $\mD_1$ and $\mD_2$,
intersecting at  a node,  denoted by  $q$.\newline
Consider the following two cases:
\newline \newline
Case II-a: the node $q$ is an untwisted point, i.e. $q \in \pi^{-1}(X_{reg})$.
\newline \newline
Case II-b: the node $q$ is a twisted point.
\newline \newline
Case II-a. This is easy; one of the irreducible components has no twisted point on it.
\newline \newline
We divide Case II-b into several subcases. 
It suffices to study the following subcases by symmetry:
\newline \newline
Case II-b-1: the twisted point $\infty \in \mD_1$ and $\{x_0,\;x_1\} \subset f_1(\mD_1)$.
\newline \newline
Case II-b-2: the twisted point $\infty \in \mD_2$ and 
$\{x_0,\;x_1\} \subset f_1(\mD_1)$.
\newline \newline
Case II-b-3: the twisted point $\infty \in \mD_1$,  
$x_0  \in f_1(\mD_1)$ and $x_1 \in f_2(\mD_2)$. 
\newline \newline
Case II-b-1. This subcase  is not possible by our choice 
 of the points $x_1$,  \;$x_0$ and $x_{\infty}$.
\newline \newline
Case II-b-2. 
In this case, $\mD_1$ has one twisted point and $\mD_2$ has two twisted points.
Consider  $f_1: \mD_1 \to \mX$.  
 Note that $p=f_1(q) \in \mX$,
the image of the node $q$,  is a stacky point on $\mX$.  

Since the image  of $\mD_1$ contains $x_1$ and $x_0$ and the class $\beta$ 
is unbreakable, 
 there are only finitely many such twisted 
stable maps thanks to  Lemma ~\ref{moribandb}.  
Since every  such twisted stable map has only finitely many twisted points (the image  only intersects the 
stacky locus at finitely many points), it 
 follows that there are only finitely many such  stacky points $p$ (since  $p$
 is  the image of a twisted point on the source curve). 

Recall  that $\mD_2$ has two stacky points, denoted by $\{q,\; \infty\}$, such that $(\pi \circ f_2)(q)=\pi(p)$ and
 $(\pi \circ f_2)(\infty)= \pi(x_{\infty})$.  
Consider the morphism  $f_2: \mD_2 \to \mX$ as an element of $\mK_{0,2}(\mX, \beta)$. 
Consider  the set of all possible twisted stable maps $[f] \in \mK_{0,2}(\mX, \beta)$ such 
 that $\pi \circ f(q)= \pi(p)$ and $\pi \circ f(\infty)= \pi (x_{\infty})$. 
This set is a finite set by Lemma ~\ref{moribandb2}.  
Since (1) $p$ (the image of the node)  is chosen from a finite set,  and (2)  for every   $p$  there are only finitely 
many such twisted stable maps $[f] \in \mK_{0,2}(\mX, \beta)$, it follows that there are only finitely many such twisted stable maps $f_2: \mD_2 \to \mX$.  
This  shows that  Case II-b-2 also forms a finite set. 
\newline \newline
Case II-b-3.  

The image of $f_1: \mD_1 \to \mX$ contains $x_{\infty}$ and $x_0$.
There are only finitely many such $f_1:D_1 \to \mX$ by 
Lemma ~\ref{moribandb}. It follows that the set of all possible nodes $p$ is also 
finite as  in Case II-b-2. 

The  image of $f_2:\mD_2 \to \mX$ contains $p$ and $x_1$.
By Lemma ~\ref{moribandb} again, the set of all such $f_2: \mD_2 \to \mX$ is also finite. 
\newline \newline
Step III. Concluding the proof:

By Step II,   all  bad cases, i.e.  Case I-a, Case I-b, Case II-b-1, Case II-b-2 and  Case II-b-3, form a finite set. 
Denote  this finite set by $S$. 
Since $dim \mK \geq 2$, we can find a proper irreducible curve $T \to \mK$ which is finite to its image and does not meet the finite 
 set $S$.  Pull back the  family over $\mK$ to $T$. 
It follows that  Case II-a is the only possible type of degeneration  in this family.  This concludes the proof.     
\end{proof}

The next two lemmas are needed in the proof of Proposition ~\ref{general}.
\begin{Lem}\label{moribandb}
Let $\beta \in H_2(X)$ be a homology class such 
that $\beta$ is unbreakable, i.e. it can not be written as the 
sum of at least two (effective) curve classes.  Let $x_1$ and $x_2$ be two distinct points on $\mX$  and 
$x_1 \in \pi^{-1}(X_{reg})$. 
Let $k=1,\;2$.
Then there are only finitely many  stable maps $[h] \in \mK_{0,k}(\mX, \beta)$ such that $x_1 \in h(\mC)$ and $x_2 \in h(\mC)$.
\end{Lem}
\begin{proof}
The proofs of $k=1$ case and $k=2$ case are similar; we only prove $k=1$ case here.
For any $[f] \in \mK_{0,0}(X, \beta)$, the domain curve is irreducible since $\beta$  is unbreakable.  
By bend and break, there are only finitely many stable maps $[f_i]\in \mK_{0,0}(X, \beta)$ whose image 
contains $\pi(x_1)$ and $\pi(x_2)$. 
Denote this finite collection  by $I=\{ [f_i]| i \in I \}$. 
Since $x_1 \in \pi^{-1}(X_{reg})$,  the image $f_i(\PP^1)$  can only intersect $X_{sing}$ at finitely 
many points. 

Consider the composition morphism $$ \mK_{0,1}(\mX, \beta) \to \mK_{0,1}(X, \beta) \to \mK_{0,0}(X, \beta),$$
 where  the  morphism $\mK_{0,1}(\mX, \beta) \to \mK_{0,1}(X, \beta)$ is  quasi-finite \cite{av02}, and the morphism $\mK_{0,1}(X, \beta) \to \mK_{0,0}(X, \beta)$ (forgetting the marked  point) is projective. 
Consider any twisted stable map $[f] \in \mK_{0,1}(\mX, \beta)$. Let $q$ be the stacky point on the 
domain curve.  Since a stacky point can only be mapped to a stacky point,  the image $p:=f(q)$ lies in the stacky locus of $\mX$ and $\pi(p) \in X_{sing}$.
Note that one can only endow non-trivial stack structures on these points. 
 Therefore, for every $[f_i] \in I$, there are  only finitely many ways to 
endow stack structures on the source curve 
such that we can lift the stable map  (into $X$) to a twisted stable map (into  $\mX$). This concludes the proof.

\end{proof}

\begin{Lem}\label{moribandb2}
Let $\beta \in H_2(X)$ be the  homology class in Lemma  ~\ref{moribandb}.  Let $x_1$ and $x_2$ be 
two points on $\mX$ (not necessarily distinct).    Then there are only finitely many  stable 
maps $[h] \in \mK_{0,2}(\mX, \beta)$ such that 
 $\pi \circ h(0)= \pi(x_1)$ and $\pi \circ h(\infty)=\pi (x_2)$, where $0$ and $\infty$ are the twisted points
 on the source curve of $[h]$.
\end{Lem}
\begin{proof}
Consider the forgetful map $\mK_{0,2}(\mX, \beta) \to \mK_{0,2}(X, \beta)$.
Recall that this morphism is quasi-finite  and the stack 
  $\mK_{0,2}(\mX, \beta)$ is proper \cite{av02}. 
Consider the collection of all possible ($2$-pointed) stable 
maps $[f] \in  \mK_{0,2}(X, \beta)$ such 
that $f(0)=\pi ( x_1)$ and $f(\infty)= \pi (x_2)$. 
It has only finitely many elements
 by Mori's bend and break. 
This concludes the proof. 
\end{proof}
The next proposition is standard: 
\begin{Prop}\label{general2}
 Assumptions as in Assumptions ~\ref{assu}. There exists a twisted stable map $f: \PP^1 \to \mX$ such that 
its $-K_{\mX}$-degree is $n+1$ and its image does not meet $\pi^{-1} (X_{sing})$. 
\end{Prop}
\begin{proof}
By Proposition ~\ref{general}, there is a twisted rational curve $g: \PP^1 \to \mX$ whose image does not lie on $\pi^{-1}(X_{sing})$.
Pick a point $q \in g(\PP^1) \subset \mX$ which is not on 
$\pi^{-1}(X_{sing})$. 
The dimension (at the point  $[g]$) of the family of   twisted stable maps $[f] \in \mK_{0,0}(\mX,n+1)$
 whose image contains $q$ and  meets  $\pi^{-1}(X_{sing})$
 is at most $dim \; X_{sing} \leq n-2$ (by a bend and break argument), while the dimension (at $[g]$) of 
twisted stable  maps $[f] \in \mK_{0,0}(\mX, n+1)$ whose image contains  $q$ is at least $n-1$
 (by Lemma ~\ref{def} and Claim ~\ref{chow}).  
By a simple  dimension count, it is easy to see that there 
is a rational curve which does not meet $\pi^{-1}(X_{sing})$.  
\end{proof}
\begin{Rem}\label{general3}
 From Proposition ~\ref{general2}, it follows easily  that  
for a general point $x \in \mX$
there exists a curve $[g] \in \mK_{0,0}( \mX, \beta)$ whose image 
contains $x$ and does not meet  $\pi^{-1}(X_{sing})$.  
\end{Rem}

\begin{notation}\label{unbreakclass}
Let $f: \PP^1 \to \mX$ be the rational curve in Proposition ~\ref{general}.
Let $\beta= [(\pi \circ f)_*(\PP^1)] \in R \subset H_2(X)$. We will fix this
class in the rest of this section.
\end{notation}
The next lemma is a simple observation:
\begin{Lem}\label{scheme}
Assumptions and notation  as in Assumptions ~\ref{assu} and Notation ~\ref{unbreakclass}. 
The Deligne-Mumford stack $\mK_{0,0}(\mX, \beta)$ is a projective scheme. 
 \end{Lem}
\begin{proof}
Let $[f] \in \mK_{0,0}( \mX, \beta)$.  
Note that the domain curve is always irreducible  since 
 there is no marked point and
$n+1$ is the minimal possible $-K_{\mX}$-degree. 
Note that $f: \PP^1 \to \mX$ is
birational to its image in $\mX$ (again, since $n+1$ is the minimal 
possible $-K_{\mX}$-degree). This implies that there is no non-trivial
automorphism for any twisted stable map in $\mK_{0,0}( \mX, \beta)$.
It follows that $\mK_{0,0}( \mX, \beta)$ is a proper scheme. Since it is  quasi-finite
to the projective scheme $K_{0,0}( X, \beta)$. It has to be projective. 
\end{proof}

Let $[f] \in \mK_{0,0}( \mX, \beta)$  be a twisted stable map into
$\mX$ which does not meet the preimage of $X_{sing}$. Let $Z \subset
\mK_{0,0}(X , \beta)$ be an irreducible component  which contains $[f]$
and $\tilde{Z}$ the normalization. Consider the finite  morphism
$\mK_{0,0}(\mX,\beta) \to K_{0,0}(X,\beta)$. Let $Z'$ be the component
of $K_{0,0}(X,\beta)$ which contains the  image of $Z$, i.e. stable
maps  of the form $[\pi \circ h]$ with $[h] \in Z$. We also take the
normalization $\tilde{Z'}$ of $Z'$.
The next lemma compares 
these two components. It is needed in the next
section.
\begin{Lem}[\cite{ct05a} Lemma 3.12]\label{iso}
The natural map $\tilde{Z} \to \tilde{Z'}$ is an
isomorphism.
\end{Lem}

\section{Results from \cite{ke00} and \cite{ke01} and some remarks}\label{kebekus}
 
We return to   work on the variety $X$, rather than
the stack $\mX$. Most  results in this section are taken from
\cite{ke00} and \cite{ke01}. The reader can also consult 
\cite{ct05a} Section 4  for more details.   
  
Let $x$ be a general point on $X$ and $R \in \overline{NE}(X)$ a
$K_X$-negative extremal ray. Let $f:\PP^1 \to
\mX$ be a twisted rational curve such that
\begin{enumerate}
\item $\PP^1 \cdot _{f} (-K_{\mX})=n+1$,
\item $[\pi \circ f (\PP^1)] \in R$, and
\item the image $f(\PP^1)$ contains $x$ and does not
meet $\pi^{-1}(X_{sing})$.
\end{enumerate}
 The existence of such a curve follows
from Proposition~\ref{general2} and Remark ~\ref{general3}.  

Fix the homology  class  $[(\pi \circ f)_* (\PP^1)] \in R \subset H_2(X)]$.  
Let  
$[f] \in H_x \subset
\mK_{0,0}( \mX, [(\pi \circ f)_* (\PP^1)])$  be  an irreducible component of the family of
stable map through $x$. 
By choosing a suitable component (in fact, there is only one component), we may assume that 
 $$dim_{[f]} H_x  \geq \PP^1 \cdot f^*(-K_X) -dim Aut(\PP^1, 0)= n+1 - dim Aut(\PP^1, 0)=n-1.$$  

By a standard  bend and break argument (see Claim ~\ref{chow}), we have $dim H_x=n-1$. 

  Let $\tilde{H}_x \to H_x$
be the normalization. By Lemma~\ref{iso}, the variety $\tilde{H}_x$
can be viewed as an irreducible component of the subfamily of stable
maps (into $X$) passing through $x$.

Consider the diagram
$$
\begin{CD}
U_x @>{i_x}>>X\\
@V{\pi_x}VV\\
\tilde{H}_x
\end{CD}$$
where $\pi_x: U_x \to \tilde{H}_x $ is  the universal family over
$\tilde{H}_x$ and $i_x: U_x \to X$ the universal
morphism into $X$.  
Note that 
the preimage $i^{-1}(x)$ contains a section, denoted by $\sigma_{\infty}$, and at most a finite
number of other points $z_i$.

Let
 $\tilde{X}  \to X$ be the blowing up of
 $X$ along a general  point $x$. There is a rational
map $\tilde{i}_x: U_x \dashrightarrow  \tilde{X}$ lifting the
morphism $i_x: U_x \to X$.

We now sketch the proof of Theorem ~\ref{main} following  \cite{ke01}.
\newline \newline
Step 1. The variety $U_x \cong \widetilde{ \PP^n}$.   \newline 
Let $E \cong \PP(T_X^*|_x)$ be the exceptional divisor of $\tilde{X} \to X$,
the blowing up of $X$ along the point $x$.
First note that   $\tilde{i}_x |_{\sigma_{\infty}}: \sigma_{\infty} \to
\PP(T_X^*|_x)\cong \PP^{n-1}$ is a finite morphism \cite{ke00}. By
\cite{ke01} Proposition 3.1, this is indeed an isomorphism. In particular, $\tilde{H}_x$ is smooth. Since $U_x$ is a $\PP^1$-bundle over $\tilde{H}_x$, $U_x$ is also smooth. 
Note that 
the universal property of the blow up implies that $\tilde{i}_x: U_x \to X$ is a morphism.
Also by  \cite{ke01} Proposition 3.1, 
 the normal bundle $N_{ \sigma_{\infty},\;U_{x} }
 \cong N_{E,\; \tilde{X} } \cong \mO_{\PP^{n-1}}(-1)$.
It is now easy to see that $U_x \cong  \PP(\mO_{\PP^{n-1}}(-1) \oplus
\mO_{\PP^{n-1}})$ which is just  the blowing up  of $\PP^n$ along a  point.\newline 

Consider  the morphism $i_x: U_x  \to X$ and its Stein factorization:
$$\begin{CD}
 U_x  @>{\operatorname{\alpha}}>> Y @>{\operatorname{\beta}}>>X,
\end{CD}$$
where $\alpha$ has connected fibers and $\beta$ is a finite morphism. \newline
\newline
Step 2.  $Y \cong \PP^n$. \newline  
It is easy to see that $Y \cong \PP^n$ since both the morphisms  $\alpha:U_x \to Y$ and $\alpha': U_x \to \PP^n$ contract the  divisor $E \subset U_x$. 
\newline \newline 
Step 3. The finite morphism $ \beta: Y \to X$ 
is an isomorphism. 
\newline 
Write $K_Y= \beta^* K_X +R$ where $R$ is the (effective) ramification divisor.
Note that the $\beta$-image of  lines through $\alpha(x)$ are curves associated with $\tilde{H}_x$.
 For any point   $[C] \in \tilde{H}_x$ ($C$ is a rational  curve of $-K_X$-degree $n+1$),  we have 
$C \cdot (-K_X)=n+1=(line)\cdot -K_{\PP^n}$.  
This implies that the divisor $R$
 is empty (recall that any effective divisor on $\PP^n$ is ample) and  $K_Y= \beta^* K_X$.

Suppose  the degree of the morphism $\beta$ is $d$. 
Take a general line $l \subset \PP^n \cong Y$. 
Let $C' \subset X$ be  the image of $l$ (under the morphism $\beta$).
The class $\beta_* [l]=d [C']$. Since  $-K_Y= \beta^* (-K_X)$, it follows that $$n+1=l \cdot (-K_Y)= l \cdot \beta^*(-K_X)=dC' \cdot (-K_X) \geq d(n+1).$$   This shows that $d=1$ 
 and $\beta$ is  an isomorphism.

\begin{Rem}
 There is a shorter, but less elementary, proof.  It  uses  the main theorem in \cite{cmsb02}. 
Recall the main theorem from
 \cite{cmsb02}:
\begin{Thm}[=\cite{cmsb02} Main Theorem 0.1]\label{cmsbmain}
Let $X$ be a normal projective variety over $\CC$. If $X$ carries a closed, maximal, unsplitting doubly dominant family of rational curves, then $X$ is isomorphic to $\PP^n$.
\end{Thm}
Proposition ~\ref{general2}   implies the existence of such a family of rational curves.

\end{Rem}
  \begin{Rem}
If $X$ has at worst LCIQ singularities, we can prove the same result with an
extra assumption that $dim\; X_{sing} < (n-1)/2$. H.-H. Tseng and I are working
on removing this extra assumption.
\end{Rem}


\begin{thebibliography}{HLOY02}
\bibitem[ACV03]{acv} D. Abramovich, A. Corti, and A. Vistoli, \emph{Twisted bundles and admissible covers}, Comm.\ Algebra 31 (2003) 3547--3618.

\bibitem[AV02]{av02} D. Abramovich, A. Vistoli. \emph{Compactifying the space of stable maps}, J. Amer. Math. Soc. 15, no. 1 27-75, 2002.

\bibitem[Am03]{am03} F. Ambro.  \emph{Quasi-log varieties}, Tr. Mat. Inst. Steklova 240 (2003), Biratsion. Geom. Linein. Sist. Konechno Porozhdennye Algebry, 220-239; translation in Proc. Steklov Inst. Math. 2003, no. 1(240), 214-233, 2003.

\bibitem[CMSB02]{cmsb02} K. Cho, Y. Miyaoka and  N.I. Shepherd-Barron.
\emph{ Characterizations of projective spaces and applications to complex symplectic manifolds}, Higher dimensional birational geometry (Kyoto, 1997),  1--88, Adv. Stud. Pure Math., 35, Math. Soc. Japan, Tokyo, 2002.
\bibitem[CT05]{ct05} J.- C. Chen and H.- H. Tseng.
 \emph{Cone Theorem via Deligne-Mumford stacks}, Preprint math.AG/0505043, 2005.\bibitem[CT05a]{ct05a} J.- C. Chen and H.- H. Tseng.
 \emph{Note on Characterization of Projective spaces}, Preprint math.AG/0509649, 2005.
\bibitem[De01]{de} O. Debarre.   \emph{Higher-dimensional algebraic geometry}, Springer-Verlag, 2001.
\bibitem[Ka91]{ka91} Y. Kawamata, \emph{On the length of an extremal rational curve}, Invent. Math. 105 (1991), no. 3, 609--611.
\bibitem[Ke00]{ke00} S. Kebekus.   \emph{Families of singular rational curves},  J. Algebraic Geom. 11 (2002), no. 2, 245-256,   math.AG/0004023.

\bibitem[Ke01]{ke01} S. Kebekus.   \emph{Characterizing the projective space after Cho, Miyaoka and Shepherd-Barron}, Complex Geometry (G\"ottingen, 2000)147-155, Springer,  Berlin, 2002,  math.AG/0107069.
\bibitem[Ko96]{ko96} J. Koll\'ar. \emph{Rational curves on algebraic varieties}, Springer-Verlag, 1996.

\bibitem[KM98]{km98} J. Koll\'ar, S. Mori. \emph{Birational Geometry of Algebraic Varieties}, Cambridge Tracts in Mathematics 134, Cambridge University Press, 1998.


\end{thebibliography}
\end{document}